\documentclass[12pt]{amsart}
\usepackage{amsmath, amssymb, amsfonts}
\makeatletter
\def\ps@pprintTitle{%
   \let\@oddhead\@empty
   \let\@evenhead\@empty
   \let\@oddfoot\@empty
   \let\@evenfoot\@oddfoot
}
\makeatother

\newtheorem {theorem}{Theorem}[section]
\newtheorem {lemma}[theorem]{Lemma}

\newtheorem {remark}[theorem]{Remark}

\newcommand{\R}{\mathbb{R}}
\newcommand{\p}{\partial}
\newcommand{\ba}{\begin{array} }
\newcommand{\ea}{\end{array} }

\def\R{{\mathbb{R}}}

\newcommand{\be}{\begin }

\numberwithin{equation}{section}

\def\ba{\begin{aligned}}
\def\ea{\end{aligned}}
\def\be{\begin{equation}}
\def\ee{\end{equation}}

\def\t#1{\tilde{#1}}
\def\eps{\varepsilon}

\def\ra{\rightarrow}

\def\div{\text{div}}

\def\O{\Omega}

\def\f{\frac}
\def\p{\partial}
\def\|{\parallel}
\def\a{\alpha}
\def\g{\nabla}
\def\ld{\lambda}
\def\lap{\triangle}

\newcommand {\eqdef }{\ensuremath {\stackrel {\mathrm {def}}{=}}}

\def\u{\mathbf{u}}

\def\O{\Omega}

\def\n#1#2{\|#1\|_{#2}}
\def\D{\mathcal{D}}
\allowdisplaybreaks

\begin{document}
	\vspace{1cm}
	\title[local existence and blowup to 2d MHD]
	{Finite time blowup of strong solutions to the two dimensional MHD equations}
	\author{Xiangdi Huang, Zhouping Xin, Wei YAN$^*$}
	\thanks{* Corresponding author.}
	
	\address{Xiangdi HUANG\hfill\break\indent
		Institute of Mathematics,
		\hfill\break\indent
		Academy of Mathematics and Systems Sciences,\hfill\break\indent
		Chinese Academy of Sciences, Beijing 100190, China}
	\email{xdhuang@amss.ac.cn}
	
	\address{Zhouping XIN\hfill\break\indent
		Institute of Mathematical Sciences,
		\hfill\break\indent
		The Chinese University of Hong Kong, Shatin, N.T., Hong Kong S. A. R., China}
	\email{zpxin@ims.cuhk.edu.hk}
	 
	\address{Wei YAN\hfill\break\indent
		School of Mathematics,\hfill\break\indent
		Jilin Unversity, Changchun 130012, China}
	\email{wyanmath@jlu.edu.cn}
	
	\maketitle

\begin{abstract}
Whether the smooth solution of the multi-dimensional viscous compressible fluids will blow-up in finite time has always been a chanllenging problem. In the recent work\cite{FM}, Merle et al.  proved that there are smooth solutions to the 2D radially symmetric compressible Navier-Stokes equations which will inevitably form shell singularities in finite time.\\
\indent In this article, we first prove the existence of local strong solutions that allow vacuum for the two-dimensional viscous compressible MHD equations on bounded domains without magnetic diffusion. Furthermore,  it is shown that if the initial data are radial symmetric and its vacuum set contains a ball centered at the origin where the total magnetic field is non-trivial, then the radial symmetric strong solution to the initial boundary value problem will definitely blow up in finite time. This is the first  example for the formation of finite time singularity of strong solutions that allows interior vacuum  of a viscous compressible fluid.

\end{abstract}


\section{Introduction}

Consider the following magnetic hydrodynamic equations without magnetic diffusion,
\be\label{mhd}
\left\{\ba
&\rho_t + \div(\rho u) = 0,\\
&(\rho u)_t +\div(\rho u\otimes u) +\g p = \div(\mu\D u) + \g(\ld \div u) + (\g\times B)\times B,\quad\text{ in }\ \O\\
&B_t-\g\times(u\times B) = 0,\\
&\div B = 0.
\ea\right.
\ee
with $(0,T)\times\Omega$ with $\Omega\subset \R^N(N=2,3)$ being a bounded domain with a smooth boundary.

Here the unknown functions $\rho, u, B$ represent the density, velocity and the magnetic field respectively. The constants $\mu$ and $\ld$ are the shear viscosity and bulk viscosity respectively which satisfy the following physical condition:
\be
\mu>0, \quad\mbox{and } \ld+\f{2}{N}\mu\ge 0.
\ee
The pressure $p$ is assumed to satisfy the following $\gamma-$ law:
\be
p=a\rho^{\gamma}\quad \mbox{for some} \quad\gamma>1.
\ee
$\D u$ is the deformation tensor given by
\be
\D u=\frac{\nabla u+\nabla u^t}{2}
\ee

The system (\ref{mhd}) is supplemented with the initial conditions
\be\label{bc-1}
\rho(x,0)=\rho_0(x),u(x,0) = u_0(x),\quad B(x,0)= B_0(x)
\ee
and the boundary condition
\be\label{bc-2}
u(x,t) = 0,\quad\text{on}\ \p\O.
\ee

 There are plenty of literatures studying  strong solutions for the compressible isentropic MHD equations. For the classical model with constant viscosity coefficients and positive magentic dissipation, the local existence of strong solutions was established by  Vol’pert–Hudjaev\cite{1972Vol} without vacuum. Later, Fan-Yu\cite{2009Fan} proved the eixstence of unique local strong solutions to the 3D compressible magnetohydrodynamic equations with nonnegative thermal conductivity or infinite electric conductivity for either the whole space or bounded domains when the  initial density is allowd to vanish. The global existence of classical solutions to the Cauchy problem of the compressible isentropic MHD equations was established by Kawashima\cite{1983Ka}  as long as the initial data is a small perturbation of a non-vacuum constant. By generalizing the ideas and proofs of Huang-Li-Xin\cite{huang2012global} for the compressible Navier-Stokes equations, Li-Xu-Zhang\cite{2013lihl}  proved the global existence of classical solutions to the Cauchy problem of the compressible isentropic MHD equations in the presence of vacuum with small total energies. In the absence of magnetic resistivity,  Zhu\cite{zhu-s} established local existence of classical solutions for the Cauchy and initial boundary value problem(IBVP) to the 3D compressible isentropic MHD equations with initial vacuum where  a higher regularity of initial data is required.

On the other hand,  whether the local classical solutions can be a global one for arbitrary initial large data  remains a long-standing open problem. In their recent work, when the density is strictly positive, Merle\cite{FM}  proved that there are smooth solutions to the 2D radially symmetric compressible Navier-Stokes equations, which will inevitably form shell singularity in finite time. There are also some notable works on the blowup of classical solutions as long as vacuum is allowed initially, see \cite{2006Cho,1998Xin,2013XY} and references therein.  In particular, it should also be noted that Rozanova\cite{Ro} showed the nonglobal existence of the classical solution to  the Cauchy problem of the isentropic magnetohydrodynamic equations with finite mass and energy. However, the local existence of such a classical solution satisfying Rozanova's condition remains unknown.  It should be noted that the blowup phenomena reveald by the previsou work\cite{1998Xin, 2013XY} are always  related to the special distribution of the initial density, such as fast decay to vacuum at far field in the whole space or admitting  some isolated mass group. More precisely, the breakdown of classical solutions are caused by the dynamic behavior of external vacuum. In sharp contrast to this, there is no study on the possible blowup of classical solutions for interior vacuum.

In this manuscript, we will investigate the 2D compressible isentropic MHD equations without magnetic dissipation in bounded domains and aim to show that a class of strong solutions with nonnegative initial density will inevitably blow-up in finite time.

To achieve this goal, we first give a detailed proof of the existence of local strong solution to the IBVP for the 2D compressible isentropic MHD equations with nonnegative density.  It should be noted that the  local well-posedness theory of strong solutions to the IBVP for the 2D compressible isentropic MHD equations with nonnegative density in the  absence of magnetic resistivity has not been established yet. Though  Fan-Yu\cite{2009Fan} investigated a more complicated system and showed  the local  existence of strong solutions to the 3D compressible MHD equations with nonnegative thermal condutivity or infinite electric conductivity without magnetic resistivity for the whole space or the bounded domain with smooth boundary. Yet the isentropic  model is not a sub-system of the system in Fan-Yu\cite {2009Fan} and  the results  there do not yield directly the existence of local strong solutions to the 2D compressible isentropic MHD equations which is needed in this paper. Thus for completeness, we still present the main estimates and construct the local solutions through by a similar but a more simplified proof process as Fan-Yu\cite{2009Fan}. It is also noted that the strong solution established in our article has lower regularity in initial density compared to the classical smooth solution both for the Cauchy and IBVP problem in Zhu\cite{zhu-s}.

  With the local existence of strong solutions at hand, out next crucial step is to derive an upper bound for life span of any strong solution by analysing carefully the interactions between the vacuum and the magnetic field. The key ideas and motivations for this can be sketched as follows.
   
   The first key element is  the conservation of the total magnegic field over a vacuum ball. This can be motivated by the following elementary observations. Note that the 2-dimensional Maxwell equations for the magnetic field $\boldsymbol {\vec{B}}$ can be written as 
   \be\label{keyy}
\boldsymbol {\vec{B}}_t +\boldsymbol {\vec{u}}\cdot\nabla \boldsymbol {\vec{B}} +\boldsymbol {\vec{B}}div \boldsymbol {\vec{u}}=\boldsymbol {\vec{B}}\cdot\nabla \boldsymbol {\vec{u}}.
\ee 
Taking the inner product of  $\eqref{keyy}$ with $\frac{\boldsymbol {\vec{x}}^{\bot}}{|\boldsymbol {\vec{x}}|^2}$ leads to 
\be\label{keyy-1}
\ba
& \Big(\frac{\boldsymbol {\vec{B}}\cdot \boldsymbol {\vec{x}}^{\bot}}{|\boldsymbol {\vec{x}}|^2}\Big)_t + div\Big(\frac{\boldsymbol {\vec{B}}\cdot \boldsymbol {\vec{x}}^{\bot}}{|\boldsymbol {\vec{x}}|^2}\boldsymbol {\vec{u}}\Big) 
& =  div\Big(\frac{\boldsymbol {\vec{u}}\cdot \boldsymbol {\vec{x}}^{\bot}}{|\boldsymbol {\vec{x}}|^2}\boldsymbol {\vec{B}}\Big),
\ea
\ee
here $\boldsymbol {\vec{x}}^{\bot}=(-x_2,x_1)$ and the following fact has been used 
\be
div\Big(\frac{\boldsymbol {\vec{x}}}{|\boldsymbol {\vec{x}}|^2}\Big) = \frac{n-2}{|\boldsymbol {\vec{x}}|^2}=0\quad \mbox{for $n=2$}.
\ee

Assume that $\Sigma_t$ is a time-dependent bounded domain which is transported by an initial domain $\Sigma_0$ through the  vector field $\boldsymbol {\vec{u}}$.  We integrate \eqref{keyy-1} over $\Sigma_t$, together with the Reynolds's transportation theorem, to get
\be
\Big(\int_{\Sigma_t}\frac{\boldsymbol {\vec{B}}\cdot \boldsymbol {\vec{x}}^{\bot}}{|\boldsymbol {\vec{x}}|^2}dx\Big)_t = \int_{\partial\Sigma_t}\Big(\frac{\boldsymbol {\vec{u}}\cdot \boldsymbol {\vec{x}}^{\bot}}{|\boldsymbol {\vec{x}}|^2}\boldsymbol {\vec{B}}\cdot \boldsymbol {\vec{n}}ds\Big),
\ee
here $\boldsymbol {\vec{n}}$ is unit spatial normal of $\partial\Sigma_t$. Thus the total magnetic field over $\Sigma_t$  is conserved in the sense
\be\label{bdcc}
\int_{\Sigma_t}\frac{\boldsymbol {\vec{B}}\cdot \boldsymbol {\vec{x}}}{|\boldsymbol {\vec{x}}|^2}dx =\int_{\Sigma_t}\frac{\boldsymbol {\vec{B_0}}\cdot \boldsymbol {\vec{x}}}{|\boldsymbol {\vec{x}}|^2}dx.
\ee
provided that
 
\be\label{bdc-0}
\frac{\boldsymbol {\vec{u}}\cdot \boldsymbol {\vec{x}}^{\bot}}{|\boldsymbol {\vec{x}}|^2}\boldsymbol {\vec{B}}\cdot \boldsymbol {\vec{n}} =0\quad \mbox{on $\partial\Sigma_t$.}
\ee

Note that, in particular, if the solution is radial symmetric, then \eqref{bdc-0} holds trivially for $\Sigma_t$ being a disc centered at the origin since  ($\boldsymbol {\vec{u}},\boldsymbol {\vec{B}}$) have the form 
\be\label{rasy-1}
\boldsymbol {\vec{u}}(\boldsymbol {\vec{x}},t)=\frac{\boldsymbol {\vec{x}}}{r}u(r,t), \boldsymbol {\vec{B}}(\boldsymbol {\vec{x}},t)=\frac{\boldsymbol {\vec{x}}^{\bot}}{r}B(r,t), \quad r=|\boldsymbol {\vec{x}}|.
\ee
Consequently, \eqref{bdcc} yields the following conservation law: 
\be\label{bdcc-1}
\int_{\Sigma_t}\frac{B(r,t)}{r}dx =\int_{\Sigma_0}\frac{B_0(r)}{r}dx
\ee
which will be rigirously established in Section 3 and $\Sigma_t$ will be choosen to be the evolution of an initial vacuum disc centered at the origin. The next key element in our analysis is based on the detailed balance between the viscous stress and the Lorentz force due to the magnetic field. Indeed, for a radial symmetric solution in the vacuum region, the momentum equations become 
\be\label{ub-1}
(2\mu+\lambda)\p_r \div \boldsymbol {\vec{u}} = B(B_r+\frac{B}{r}),\quad\text{for any }r\in [0, R(t)] \quad a.e.
\ee

The desired upper bound estimate on the life span of the solution can be derived by using the fact that the Lorentz force on the right hand side of \eqref{ub-1} is due to magnetic field which admits a conserved quantity as \eqref{bdcc-1}, while the left hand side of \eqref{ub-1} is the viscous stress which has suitable time decay due to the basic energy inequality. To achieve this, we will use a fractional moment argument for the momentum equation \eqref{ub-1} , which can be sketched as follows. We will look for a suitable multiplier  $f(r)$ for \eqref{ub-1} with $f(0)=f(R)=0, R=R(t)$. Then \eqref{ub-1} implies
\be\label{ff-1}
-\int_0^R (2\mu+\lambda)f'(r)div\boldsymbol {\vec{u}} dr = \int_0^R f(r)\Big(\p_r\frac{B^2}{2}+\frac{B^2}{r}\Big)dr.
\ee

Note that the left hand side of \eqref{ff-1} is bounded from above by
\be\label{ff-11}
C_1 \|div \boldsymbol {\vec{u}}\|_2\Big(\int_0^R|f'(r)|^2r^{-1}\Big)^{\frac{1}{2}},
 \ee
while the right hand side of \eqref{ff-1} is bounded from below by
\be\label{ff-12}
C_2\Big(\int_0^RBdr\Big)^2\Big(\int_0^R\Big[\frac{f(r)}{r}-\frac{f'(r)}{2}\Big]^{-1}dr\Big)^{-1}
\ee
with $C_1$ and $C_2$ being uniform positive constants. It follows from \eqref{ff-1}-\eqref{ff-12} and integrating in time that an upper bound for the life span of the strong solution can be obtained as long as one can choose the multiplier $f(r)$ such that

\be\label{ff-13}
\Big(\int_0^R|f'(r)|^{2}r^{-1}\Big)^{\frac{1}{2}}\Big(\int_0^R\Big[\frac{f(r)}{r}-\frac{f'(r)}{2}\Big]^{-1}dr\Big)\lesssim 1.
\ee
This will be achieved by choosing
  \be
  f(r)=(R-r)r^{\alpha}
  \ee
for any given constant $\alpha\in (1,2)$ in Section 3. Thus we can derive a finite upper bound for the life span of the strong solution.

This yields the first result on the formation of finite time singularities of strong solutions that allows interior vacuum for a viscous compressible fluid.
\section{Local existence}

We first prove that initial boundary value problem \eqref{mhd},\eqref{bc-1},\eqref{bc-2} has a unique local strong solution in the presence of vacuum when $\Omega\subset \R^2$ is a 2D bounded domain with smooth boundary as follows.

\begin{theorem}\label{local} ({\bf Local existence of strong solution})
	Assume that the initial data $(\rho_0, u_0, B_0)$ satisfy the following conditions
	\be\left\{\ba
	&0\le\rho_0\in W^{1,q}, B_0\in W^{1,q}, \ div B_0=0 \text{ in } \Omega,\\
	&u_0\in H^2(\O),\ u_0\big|_{\p\O} = 0,
	\ea\right.\ee
	for some constant $q>2$, and 
	\be
	\g p(\rho_0) - \div(\mu\D( u_0)) - \g(\ld \div u_0) + B_0\times(\g\times B_0) = \rho_0^{\f 12}g,\text{ in }\O,\text{ for some }g\in L^2.
	\ee
	Then, there exist a small time $T^* >0$ and a unique strong solution $(\rho, u, B)$ to the initial boundary problem \eqref{mhd},\eqref{bc-1},\eqref{bc-2} such that
	\be\left\{\ba
	& (\rho, B)\in C([0, T^*]; W^{1,q}),\ u\in C([0, T^*]; H^2)\cap L^2(0, T^*; W^{2,q}),\\
	& (\rho_t, B_t)\in C([0, T^*]; L^q), u_t \in L^2(0, T^*; H^1),\\
	& \sqrt{\rho} u_t\in L^{\infty}(0, T^*; L^2).
	\ea\right.\ee
\end{theorem}
\begin{remark}
 The proof of Theorem \ref{local} will be postponed to the last part.
\end{remark}
\begin{remark}
In the next section, we prove that, in general, the strong solution of \eqref{mhd} guaranteed by Theorem \ref{local} will breakdown in finite time in the presence of interior vacuum region.
\end{remark}

To this end, we also give the local existence of  classical solutions as long as the initial data enjoy more regularity. This was first proved by Zhu\cite{zhu-s} for the 3D case, where the domain can be a bounded domain in $R^3$ with smooth boundary or the whole space $R^3$, which can be generalized to  2D bounded domains as follows:

\begin{theorem}\label{local-c} ({\bf Local existence of classical solution})
	Assume that $q\in (2,\infty)$ is any given constant, and supplement the system \eqref{mhd}  with the following initial-boundary conditions
	\be\label{bc-3}
	(B,\rho,u)|_{t=0}=(B_0(x),\rho_0(x),u_0(x)),\quad \text{$x\in\Omega$},\quad u\big|_{\p\O} = 0.
	\ee
	
	Assume also that the initial data $(B_0, \rho_0, u_0 ,P_0)$ satisfy the following regularity conditions
	\be\label{bc-4}
	\left\{\ba
	&(B_0, \rho_0,u_0,P_0)\in W^{2,q}, \rho_0\ge 0, \\
	&u_0\in H^2,\ u_0\big|_{\p\O} = 0,
	\ea\right.\ee
	and the following compatibility condition
	\be
	\g P_0 - \div(\mu\D( u_0)) - \g(\ld \div u_0) + B_0\times(\g\times B_0) = \sqrt{\rho_0}g,\text{ in }\O,\text{ for some }g\in L^2.
	\ee
	Then, there exist a small time $T^* >0$ and a unique classical solution $(\rho, u, B)$ to the initial boundary problem \eqref{mhd},  \eqref{bc-3}, \eqref{bc-4} such that
	\be\label{2bc}
	\left\{\ba
	& (B,\rho, P)\in C([0, T^*]; W^{2,q}),\\
	& u\in C([0, T^*]; H^2)\cap L^2(0, T^*; H^3),\ u_t\in L^2(0,T^*; H_0^1)\\
	& tu\in L^{\infty}([0, T^*]; W^{3,q}), u_t \in L^2(0, T^*; H^1),\\
	&tu_t\in L^{\infty}(0, T^*; H^2),\ tu_{tt}\in L^{\infty}(0, T^*; H_0^1)
	\ea\right.\ee
\end{theorem}
\begin{remark}
The solution obtained by Theorem \ref{local-c} is indeed classical in $\Omega\times(0,T^*]$. It suffices to prove the following regularity for any $0<\tau<T^*\le\infty$
\be
(\rho_t,u_t,\nabla^2u)\in C(\Omega\times[\tau, T^*]).
\ee
First, by the regularity \eqref{2bc} and the emdedding Theorem, it holds that
\be
(\nabla u, \nabla^2u)\in C([\tau, T^*];L^2)\cap L^\infty([\tau, T^*];W^{1,q})\hookrightarrow C(\Omega\times[\tau, T^*])
\ee 
Then the continuity of the $\rho_t$ follows immediately from the mass equations
\be
\rho_t +u\cdot\nabla\rho+\rho\div u=0.
\ee

Similarly,  \eqref{2bc} implies
\be
u_t\in L^\infty([\tau, T^*]; H^2),\ u_{tt}\in L^2(\tau, T^*; H^1)
\ee
Consequently, for any $0<\tau\le t_1<t_2\le T^*$, one has
\be
\ba
\|\nabla u_t(t_2)-\nabla u_t(t_1)\|_{L^2}&=\|\int_{t_1}^{t_2}\nabla u_{tt}(s)ds\|_{L^2}\\
&\le \int_{t_1}^{t_2}\|\nabla u_{tt}(s)\|_{L^2}ds\\
&\le \big( \int_{t_1}^{t_2}\|\nabla u_{tt}(s)\|^2_{L^2}ds \big)^{\frac{1}{2}}(t_2-t_1)^{\frac{1}{2}}\\
&\le C(t_2-t_1)^{\frac{1}{2}}
\ea
\ee
Hence,
\be
u_t\in C([\tau,T^*];H^1)\cap L^\infty([\tau, T^*];H^2)\hookrightarrow C(\Omega\times[\tau, T^*]).
\ee
\end{remark}

\begin{remark}
Compared to the strong solution established by Theorem \ref{local}, the higher regularities in Thoerem \ref{local-c} require $P_0\in W^{2,q}$, which is equivalent to
\be
\p_{i,j}P_0=a\gamma\rho_0^{\gamma-1}\p_{ij}\rho_0+a\gamma(\gamma-1)\rho_0^{\gamma-2}\p_i\rho_0\p_j\rho_0\in L^q.
\ee
This holds automatically for $\gamma\ge 2$ but imposes extra restrictions on the intial density when $\gamma\in (1,2)$, such as 
\be
\rho_0^{\gamma-2}\p_i\rho_0\p_j\rho_0\in L^q.
\ee
For instance,  it suffices to impose
\be
\nabla \rho_0^{\gamma-1}\in L^q\quad \mbox{or}\quad \rho_0^{\gamma-2}\in L^q,\quad \mbox{$\gamma\in (1,2)$}.
\ee

\end{remark}

\section{Finite time blow-up}

In this section, we prove that the strong solution to \eqref{mhd} on bounded domains, in general, blows up in finite time in the presence of interior vacuum and non-trivial magnetic field. To this end, we consider the radially symmetric solutions of \eqref{mhd}. In this case,  we look for a class of radially symmetric functions ({\boldmath ${\rho}, \vec{u}, \vec{B}$})  with the following form
\be\label{rasy}
\boldsymbol {\rho}(x,t)=\rho(r,t), \boldsymbol {\vec{u}}(x,t)=\frac{x}{r}u(r,t), \boldsymbol {\vec{B}}(x,t)=\frac{(-x_2,x_1)}{r}B(r,t).
\ee
Here $(\rho, u, B)$ are scalar functions of $(r=|x|,t)$,$x=(x_1,x_2)$. 

Hence, in a disc $\Omega=B_{R_0}$ centered at the origin with radius $R_0$, the original system \eqref{mhd} can be rewritten as

\be\label{mhd-r}
\left\{\ba
&\rho_t + (\rho u)_r + \frac{\rho u}{r}= 0,\\
&(\rho u)_t + (\rho u^2)_r + \frac{\rho u^2}{r} +  p_r =  (2\mu+\ld)(u_r+\frac{u}{r})_r - B(B_r+\frac{B}{r}),\\
&B_t + (uB)_r = 0.
\ea\right.
\ee
with $0\leq r\leq R_0<+\infty$ and the initial data
\begin{equation}\label{initial data2}
	\rho(r,0)=\rho_{0}(r), u(r,0)=u_{0}(r), B(r,0)=B_0(r),~~~~0\leq r\leq R_0,
\end{equation}
and the Dirichlet boundary condition
\begin{equation}\label{BD2}
	u(R_0,t)=0,~~~~~t\geq 0.
\end{equation}
Note also that the continuity of the velocity and magentic field in the center will force 
\be\label{BD3}
u(0,t)=B(0,t)=0,~~~~~t\geq 0.
\ee

The main conclusion of this paper claims that the solution to the problem  \eqref{mhd-r}-\eqref{BD3} will develop finite time singularity for a class of initial data which admits an interior region of vacuum containing the origin where the magnetic field is non-trivial.

\begin{theorem}\label{blowup} ({\bf Blowup of strong solutions})
	Assume that $\Omega=B_{R_0}$ is a 2D disc centered at the origin with radius $R_0$. Assume also that the initial data ({\boldmath ${\rho_0}, \vec{u_0}, \vec{B_0}$}) of the initial boundary problem \eqref{mhd},\eqref{bc-1},\eqref{bc-2} satisfy the following regularity condition
	\be\left\{\ba
	&\boldsymbol {({\rho_0}, \vec{B_0})}\in W^{1,q},\ \boldsymbol {{\rho_0}}\ge 0,\ div \boldsymbol {\vec{B_0}}=0 \text{ in } \Omega=B_{R_0},\\
	&\boldsymbol {\vec{u_0}}\in H^2(\O), \boldsymbol {\vec{u_0}}\big|_{\p\O} = 0,
	\ea\right.\ee
	for some constant $q>2$, and 
	\be
	\g p(\boldsymbol {{\rho_0}}) - \div(\mu\D( \boldsymbol {\vec{u_0}})) - \g(\ld \div \boldsymbol {\vec{u_0}}) + \boldsymbol {\vec{B_0}}\times(\g\times \boldsymbol {\vec{B_0}}) = \sqrt{\boldsymbol {{\rho_0}}}\boldsymbol {\vec{g}},\text{ in }\O,\text{ for some }\boldsymbol {\vec{g}}\in L^2.
	\ee
	Then by Theorem \ref{local}, the initial boundary problem \eqref{mhd},\eqref{bc-1},\eqref{bc-2}  will admit a unique local strong solution ({\boldmath ${\rho}, \vec{u}, \vec{B}$}) for some time $T>0$.
	
	If, in addition that
	\be\label{ra-1}
	\boldsymbol {{\rho_0}}(x)=\rho_0(r), \boldsymbol {\vec{u_0}}(x)=\frac{x}{r}{u_0}(r), \boldsymbol {\vec{{B}}_0}(x)=\frac{(-x_2,x_1)}{r}B_0(r).
	\ee
	The problem \eqref{mhd},\eqref{bc-1} and \eqref{bc-2} will have a unique radially symmetric strong solutions $(\rho,u,B)$ for some time $T$ which satisfy the following
	\be\label{ra-2}
	\boldsymbol {{\rho}}(x,t)=\rho(r,t), \boldsymbol {\vec{u}}(x,t)=\frac{x}{r}u(r,t), \boldsymbol {\vec{B}}(x,t)=\frac{(-x_2,x_1)}{r}B(r,t),
	\ee
	where $(\rho, u, B)$ are scalar functions.
	
	Furthermore, if there exists a constant $r_0\in (0,R_0)$ such that
	\be
	\rho_0=0\quad\mbox{in $B_{r_0}$}
	\ee
	and 
	\be\label{b-1}
	\int_0^{r_0} B_0(r) dr  \neq 0.
     \ee
	Then, the local strong solution ({\boldmath ${\rho}, \vec{u}, \vec{B}$}) to the initial boundary problem \eqref{mhd},\eqref{bc-1},\eqref{bc-2} will breakdown in finite time. 
	
	More precisely, for any given constant $\alpha\in (1,2)$, the life span of the strong solution will be bounded as follows
\be
T\lesssim  \left(\frac{1}{(2\mu+\lambda)R_0}\f{(2-\a)^2}{2\Big[\f{\a}{\sqrt{2\a-2}}+  \f{\a+1}{\sqrt{2\a}} \Big]}\right)^{-2}\left(\int_0^{r_0} B_0(r) dr\right)^{-4}E_0,
\ee
here $E_0$ is the initial energy given by
\be
E_0=\frac{1}{2}\int(\boldsymbol {{\rho}_0}|\boldsymbol {\vec{u}_0}|^2+\boldsymbol {\vec{B}_0}^2)dx+\int\frac{a\boldsymbol {{\rho}_0}^\gamma}{\gamma-1}dx.
\ee
	
\end{theorem}
\begin{remark}
Note that there are two crucial assumptions in Theorem \ref{blowup}: one is that the initial density contains an open interior region of vacuum, and the other is total mangetic field is nontrivial in the vacuum region. Thus the finite blowup is a consequence of the interaction between non-zero magnetic fields and internal vacuum.To our best knowledge, this is the first study on the influence of the dynamics of internal vacuum on the evolution of strong solutions.
\end{remark}

\begin{remark}
It is quite surprising that in the absence of the magnetic field, i.e, $B\equiv 0$, it remains open whether the conclusions of Theorem 1.2  hold, which is different from the case that the initial density is of isolated mass group \cite{2013XY}.
\end{remark}

\begin{remark}
The results also hold for the classical solutions of system \eqref{mhd} guranteed by Theorem \eqref{local-c}.
\end{remark}

\proof  The proof will be divided into three steps. 

{\it Step 1.}  Momentum balance in vacuum region.

By assumption,  ({\boldmath ${\rho}, \vec{u}, \vec{B}$}) is the local radially symmetric strong solution to \eqref{mhd},\eqref{ra-1} and \eqref{ra-2} for some time $T>0$ which admits the following regularities
\be
 (\boldsymbol {\rho}, \boldsymbol {\vec{B}})\in C([0, T]; W^{1,q}),\ \boldsymbol {\vec{u}}\in C([0, T]; H^2)\cap L^2(0, T; W^{2,q}).
 \ee
 Thus as $\nabla\boldsymbol {\vec{u}}\in L^2(0, T; W^{1,q})\subset L^2(0,T;L^\infty)$, these ensure the existence and uniqueness of particle trajectories $\boldsymbol {\vec{X}}(x,t)$ satisfying
 \be
 \left\{\ba
 &\frac{\p \boldsymbol {\vec{X}}}{\p t}(x,t)=\boldsymbol {\vec{u}}(X(x,t),t),\\
 &\boldsymbol {\vec{X}}(x,0)=x
 \ea\right.
 \ee
Due to the radially symmetry, there exists a scalar function $X(r,t)$ such that
\be
\boldsymbol {\vec{X}}(x,t)=X(r,t)\frac{x}{r}.
\ee
Now let $R(t)$ be a particle path at time t in the region of interior vacuum starting from point $r_0$, that is, 
\be
\rho=0 \quad\mbox{on} \quad[0, R(t)],\quad R(0)=r_0.
\ee

 Then, in the vacuum region $[0,R(t)]$, the momentum equation in \eqref{mhd-r} becomes
\be\label{key-11}
(2\mu+\lambda)(u_r+\frac{u}{r})_r = (2\mu+\lambda)\p_r \div \boldsymbol {\vec{u}} = B(B_r+\frac{B}{r}),\quad\text{for any }r\in [0, R(t)] \quad a.e.
\ee
In the following, we just use $R$ to denote $R(t)$ for simplicity. 

{\it Step 2.} We will use a fractional moment argument for  \eqref{key-11} to detect the detailed balance between the viscous stress and the magnetic Lorentz force to derive the following inequality in a vacuum disc,  which is one of the key steps in deriving the upper bound of the life span of strong solutions. Actually, for any given constant $\alpha\in (1,2)$, we claim that
\be\label{key-c0}
\f {2-\a} 2\int_0^R B^2R r^{\a-1} dr \leq (2\mu+\lambda)\Big[\f{\a}{\sqrt{2\a-2}}+  \f{\a+1}{\sqrt{2\a}} \Big]R^{\a}\|\div \boldsymbol {\vec{u}}\|_2.
\ee
As discussed in the introduction, the main idea  to prove the above claim is to find a suitable multiplier for \eqref{key-11}. This turned out to be one of the main observations in this article.

Indeed, for any $1<\a<2$, multiplying \eqref{key-11} by $(R r^\a - r^{\a+1}) $ gives
\[
(2\mu+\lambda)(R r^\a - r^{\a+1}) \p_r\div \boldsymbol {\vec{u}}  = (R r^\a - r^{\a+1}) (\p_r\f{B^2} 2 + \f {B^2}r)\quad \mbox{a.e.}
\]
Integrating the above equation over $[0, R(t)]$ yields
\be\label{key-1}
\ba
(2\mu+\lambda)\int_0^R (R r^\a - r^{\a+1}) \p_r\div \boldsymbol {\vec{u}} dr = \int_0^R  (R r^\a - r^{\a+1}) (\p_r\f{B^2} 2 + \f {B^2}r) dr.
\ea
\ee

Note that the left hand side (LHS) of \eqref{key-1} is
\be\label{key-2}
\ba
LHS &=(2\mu+\lambda) \int_0^R (R r^\a - r^{\a+1}) \p_r\div \boldsymbol {\vec{u}} dr \\
&= -(2\mu+\lambda)\int_0^R (\a R r^{\a-1} - (\a+1)r^{\a}) \div \boldsymbol {\vec{u}} dr.
\ea
\ee

In order to show that \eqref{key-2} holds for the strong solution $(\boldsymbol {\vec{u}, \vec{B}})$, one needs to verify the regularities of the integrand. 

Indeed, recalling $\boldsymbol {\vec{u}}\in C([0,T];H^2)$ and
\be
\p_r\div \boldsymbol {\vec{u}} =\nabla\div \boldsymbol {\vec{u}}\cdot\frac{x}{r}
\ee
one has immediately that
\be
\|r^{\frac{1}{2}}\p_r\div \boldsymbol {\vec{u}}\|_{L^2(0,R)}\lesssim \|\nabla\div \boldsymbol {\vec{u}}\|_{L^2(\Omega)}\le C
\ee
The integrand function of the left hand of \eqref{key-2} can be rewritten as 
\be
(R r^\a - r^{\a+1}) \p_r\div \boldsymbol {\vec{u}} = \Big(R r^{\a-\frac{1}{2}} - r^{\a+\frac{1}{2}}\Big) \Big(r^{\frac{1}{2}}\p_r\div \boldsymbol {\vec{u}}\Big) \in L^\infty((0,T);L^2(0,R)).
\ee
On the other hand, the fact $\boldsymbol {\vec{u}}\in C([0,T];H^2)$ implies that $div \boldsymbol {\vec{u}}\in C([0,T];H^1)\subset C([0,T];L^p)$ for any $p\in (2,\infty)$. Thus it holds that
\be
\|r^{\frac{1}{p}}\div \boldsymbol {\vec{u}} \|_{L^p(0,R)}\lesssim \|\div \boldsymbol {\vec{u}}\|_{L^p(\Omega)}.
\ee

Consequently, taking p large engough that $\frac{1}{p}<\a-1$, the integrand of the righthand side of \eqref{key-2} can be rewritten as 
\be\ba
&\Big(\a R r^{\a-1} - (\a+1)r^{\a}\Big) \div \boldsymbol {\vec{u}} \\
&=\Big( \a Rr^{\a-1-\frac{1}{p}}- (\a+1)r^{\a-\frac{1}{p}}\Big) \Big(r^{\frac{1}{p}}\div \boldsymbol {\vec{u}}\Big)\in L^\infty((0,T);L^p(0,R)).
\ea
\ee
Thus \eqref{key-2} can be easily verified by an approximating sequences of smooth solutions. Indeed, say for any $\epsilon>0$, one can decompose the integral in to two parts
\be\label{epsilon}
\int_0^R = \int_0^{\epsilon} + \int_{\epsilon}^R,
\ee
then \eqref{key-2} can be derived by  $\epsilon\rightarrow 0$.

Therefore, the left hand side of \eqref{key-1} can be bounded by
\be\label{lt-1}
\ba
\frac{LHS}{2\mu+\lambda}&\leq \a R\int_0^R  r^{\a-1}|\div \boldsymbol {\vec{u}} |dx +  (\a+1)\int_0^R r^{\a} |\div \boldsymbol {\vec{u}}| dr\\
&\leq \a R\|\div \boldsymbol {\vec{u}}\|_2(\int_0^R  |r^{\a-2}|^2 rdr)^{\f 12} +  (\a+1)\|\div \boldsymbol {\vec{u}}\|_2 (\int_0^R |r^{\a-1}|^2 r dr)^{\f 12}\\
&=\a R\|\div \boldsymbol {\vec{u}}\|_2\f{1}{\sqrt{2\a-2}}R^{\a-1} +  (\a+1)\|\div \boldsymbol {\vec{u}}\|_2\f{1}{\sqrt{2\a}} R^{\a}\\
&=\Big[\f{\a}{\sqrt{2\a-2}}+  \f{\a+1}{\sqrt{2\a}} \Big]R^{\a}\|\div \boldsymbol {\vec{u}}\|_2
\ea
\ee
Similarly, the righthand side(RHS) of \eqref{key-1} admits the following estimates:
\be\label{rt-1}
\ba
RHS &= \int_0^R  (R r^\a - r^{\a+1}) (\p_r\f{B^2} 2 + \f {B^2}r) dr\\
&=-\int_0^R\f{B^2} 2 (R\a r^{\a-1} - (\a+1)r^\a)dr + \int_0^R B^2 (R r^{\a-1} - r^{\a})dr \\
&=\int_0^R (1-\f \a 2)B^2R r^{\a-1} dr + \int_0^R (\f{\a+1}{2}-1)B^2 r^{\a}dr\\
&=\f {2-\a} 2\int_0^R B^2R r^{\a-1} dr + \f{\a-1}{2}\int_0^R B^2 r^{\a}dr\\
&\geq \f {2-\a} 2\int_0^R B^2R r^{\a-1} dr.
\ea
\ee
The identity in \eqref{rt-1} holds also for the strong solution $\boldsymbol {\vec{B}}$ by applying a same approximation procedures as for the lefthand side of \eqref{key-1}. The key step is to show the following regularity
\be
(R r^\a - r^{\a+1}) \p_r\f{B^2} 2 = \Big(R r^{\a-\frac{1}{q}} - r^{\a+1-\frac{1}{q}}\Big) \Big(r^{\frac{1}{q}}\p_r\f{B^2} 2\Big)\in L^\infty(0,T;L^q(0,R)),
\ee
where one has used $\boldsymbol {\vec{B}}\in C([0,T];W^{1,q})\subset C([0,T]\times\Omega)$ with
\be
B=\boldsymbol {\vec{B}}\cdot\frac{x^{\perp}}{r}.
\ee
Consequently,
\be
\|r^{\frac{1}{q}} \p_r\f{B^2} 2\|_{L^q(0,R)}\lesssim \|\boldsymbol {\vec{B}}\|_{L^\infty(\Omega)}\|\nabla \boldsymbol {\vec{B}}\|_{L^q(\Omega)}\le C.
\ee
Meanwhile, 
\be
B^2r^{\a-1} + B^2r^{\a}\lesssim C\in L^{\infty}((0,T)\times\Omega).
\ee

Combing \eqref{lt-1} with \eqref{rt-1} yields
\be\label{key-c}
\f {2-\a} 2\int_0^R B^2R r^{\a-1} dr \leq (2\mu+\lambda)\Big[\f{\a}{\sqrt{2\a-2}}+  \f{\a+1}{\sqrt{2\a}} \Big]R^{\a}\|\div \boldsymbol {\vec{u}}\|_2.
\ee

{\it Step 3.} Conservaton of total magnetic field in vacuum region.

The next key observation is that $B$ is conserved in the vacuum region, i.e, 
\be\label{r-10}
\Big|\int_0^{R(t)} B dr\Big| = \Big|\int_0^{R(0)=r_0} B_0 dr\Big| = C_0>0.
\ee
As long as \eqref{r-10} holds true, the left-hand side will be bounded below by some positive constant for non-trivial magnetic field satisfying \eqref{b-1}. This gives the lower bound of the $\|\div \boldsymbol {\vec{u}}\|_{L^2}$ as well. Consequently, the life span of the strong solution follows immediately by the basic energy inequality with disspation.

Indeed, recall the regularity of $\boldsymbol {\vec{u}}$ and $\boldsymbol {\vec{B}}$ in \eqref{2bc} that
\be
\boldsymbol {\vec{u}}\in C([0,T];H^2)\subset C([0,T]\times\Omega),\ \boldsymbol {\vec{B}}\in C([0,T];W^{1,q})\subset C([0,T]\times\Omega),
\ee
which implies
\be
u\in C([0,T]\times\Omega),\ B\in C([0,T]\times\Omega)
\ee
Integrating $\eqref{mhd-r}_3$ over $[0,R(t)]$ yields
\be\label{r-11}
\int_0^{R(t)}B_tds=\big(\int_0^{R(t)}Bds\big)_t-R'(t)B(R(t),t)
\ee
Note that $R(t)$ satisfies
\be\label{r-12}
R'(t)=u(R(t),t),\ R(0)=r_0.
\ee
Again one has
\be\label{r-13}
\int_0^{R(t)}(uB)_rdr=u(R(t),t)B(R(t),t).
\ee
Note that \eqref{r-13} can be also shown by decomposing the integral into two  parts as done in \eqref{epsilon} with the help of \eqref{BD3}.

Collecting \eqref{r-11}-\eqref{r-13} together shows
\be
\big(\int_0^{R(t)}Bdr\big)_t=0
\ee
which verifies \eqref{r-10}.

Noting $1<\a<2$, one has
\be\label{key-c1}
\ba
C_0^2 &= \Big|\int_0^{R(t)} B dr\Big|^2\leq \int_0^R B^2r^{\a-1}dr \int_0^R r^{1-a}dr\\
&=\int_0^R B^2 r^{\a-1}dr \int_0^R r^{1-a}dr\\
&=\int_0^R B^2 r^{\a-1}dr \f{1}{2-\a}R^{2-\a}.
\ea\ee

Finally, we deduce from \eqref{key-c} and \eqref{key-c1} that
\[\ba
&(2\mu+\lambda)\Big[\f{\a}{\sqrt{2\a-2}}+  \f{\a+1}{\sqrt{2\a}} \Big]R^{\a}\|\div \boldsymbol {\vec{u}}\|_2 \\
&\geq \f {2-\a} 2R\int_0^R B^2 r^{\a-1} dr\\
&\geq \f {2-\a} 2 R C_0^2(2-\a)R^{\a-2}\\
&=\f {(2-\a)^2} 2 C_0^2 R^{\a-1},
\ea\]
which implies
\[
(2\mu+\lambda)R \|\div \boldsymbol {\vec{u}}\|_2 \geq \f{C_0^2(2-\a)^2}{2\Big[\f{\a}{\sqrt{2\a-2}}+  \f{\a+1}{\sqrt{2\a}} \Big]}>0.
\]
Consequently, noting that $R(t)\leq R_0$, we have a lower bound for $\div u$ as follows
\be\label{key}
\|\div \boldsymbol {\vec{u}}\|_2 \geq \frac{1}{(2\mu+\lambda)R_0}\f{C_0^2(2-\a)^2}{2\Big[\f{\a}{\sqrt{2\a-2}}+  \f{\a+1}{\sqrt{2\a}} \Big]}>0.
\ee
Recall the basic energy estimates 
\be
\|\boldsymbol {\rho}\|^{\gamma}_{L^{\infty}L^{\gamma}}+\int_0^T\|\nabla \boldsymbol {\vec{u}}\|_{L^2}^2dt \lesssim E_0,
\ee
here $E_0$ is the initial energy.

Thus, the life span of the strong solution must satisfy
\be
T\lesssim  \left(\frac{1}{(2\mu+\lambda)R_0}\f{C_0^2(2-\a)^2}{2\Big[\f{\a}{\sqrt{2\a-2}}+  \f{\a+1}{\sqrt{2\a}} \Big]}\right)^{-2}E_0.
\ee
This finishes the proof of Theorem 1.2.
\endproof
The next sections will be devoted to the existence of local strong solutions in presence of vacuum.

\section{Linearized Problem}
Through out this section and thereafter, for simplicity, we will use $(\rho, u, B)$ instead of ({\boldmath ${\rho}, \vec{u}, \vec{B}$}) to denote the solutions to the system \eqref{mhd}.

The first step is to prove the existence of strong solution for the linearized system with some uniform estimates independent of the approximating sequences.

Consider the following linearized problem for \eqref{mhd}
\be\label{leq1}
\left\{\ba
&\rho_t + \div(\rho v) = 0,\\
&(\rho u)_t + \div(\rho v\otimes u) +\g p = \div(\mu\D(u)) + \g(\ld\div u) - B\times(\g\times B),\\
&B_t  -\g\times(v\times B) = 0.\\
&(\rho, \rho u, B)|_{t=0} = (\rho_0, \rho_0 u_0, B_0),\\
& u = 0\quad\text{ on } \p\O.
\ea\right.\ee
with $\rho_0, u_0, B_0$ satisfying
\be\label{leq2}
\left\{\ba
&2+\gamma+\n{(\rho_0, B_0)}{W^{1,q}} + \n{u_0}{H^2} +\n{g}{L^2}^2<c_0<\infty,\\
& \rho_0\geq \delta>0, \ \\
&g=\rho_0^{\f 12}(\g p(\rho_0) - \div(\mu\D( u_0)) - \g(\ld \div u_0) + B_0\times(\g\times B_0)).
\ea\right.\ee
where the known vector $v$ satisfies $v(0) = u_0,\ v\big|_{\p\O}=0$ and
\be\label{leq3}
\ba
&\n {v} {L^\infty(0, T^*; H^1)} + \beta^{-1}\n{v}{L^\infty(0, T^*; H^2)}+\n{v_t} {L^2(0,T^*; H^1)} + \n{v} {L^2(0, T^*; W^{2,q})}\leq c_1,
\ea\ee
for some fixed constant $c_0,\ c_1,\ \beta$ and $T^*$ such that
\[
1<c_0<c_1<c_2 \eqdef \beta c_1, \text{ and } 0<T^*<\infty.
\]

First, one can get the following bound of density to the linearized system. Indeed, the following results is proved in \cite{2006ChoK}. 
\begin{lemma}\label{LRHO} There exists a unique solution $\rho$ to the linearized problem such that
\be
\n{\rho}{W^{1,q}} \leq Cc_0,\ \n{\rho_t(t)}{L^q}\leq Cc_2^2.
\ee
for $0\leq t\leq \min(T^*,T_1)$ with $T_1\eqdef c_2^{-1}<1$.
\end{lemma}
\proof  As pointed by Lemma 6 in \cite{2006ChoK}, the existence of a unique solution $\rho$ was already proved by Diperna and Lions \cite{1989DL}. It following from (105) in \cite{2006ChoK} that
\[
\|\rho(t)\|_{W^{1,q}}\leq \|\rho_0\|_{W^{1,q}}\exp\Big(C\int_0^t\|\g v(s)\|_{W^{1,q}} ds\Big).
\]
Note that
\[
\int_0^t\|\g v(s)\|_{W^{1,q}} ds \leq t^{\f 12}\Big(\int_0^t \|\g v(s)\|_{H^1\cap W^{1,q}}^2 ds\Big)^{\f 12}
\leq C(c_2t + (c_2t)^{\f 12}) \leq C.
\]
for $0\leq t\leq T_1$.

We conclude that, for $0\leq t\leq T_1$, 
\[
\|\rho(t)\|_{W^{1,q}}\leq Cc_0.
\]
and
\[
\|\rho_t(t)\|_{L^q} = \|\div(\rho v)\|_{L^q}\leq C(\|\rho\|_{W^{1,q}}\|v\|_{W^{1,q}}) \leq 
Cc_0 \|v\|_{H^2} \leq Cc_0 \beta c_1\leq Cc_2^2.
\]

Moreover, for $0\leq t\leq T_1$, it holds that
\[
C^{-1}\delta \leq \rho(t,x)\leq Cc_0.
\]
\endproof

Next, we prove the existence of $B$ to the linearized system \eqref{leq1}.
\begin{lemma}\label{LB} There exists a unique solution $B$ to \eqref{leq1} such that
	\be\label{EB1}
	\|B(t)\|_{W^{1,q}}\leq Cc_0^2.
	\ee
	and
	\be\label{EB2}
	\|B_t(t)\|_{L^q}\leq Cc_0^2c_2.
	\ee
for $0\leq t\leq \min(T^*,T_1)$ with $T_1\eqdef c_2^{-1}<1$.
\end{lemma}
\proof It suffices to prove \eqref{EB1} and \eqref{EB2}. Since
\be
B_t + (v\cdot\g) B - (B\cdot\g) v + B\div u = 0,
\ee
one has, for $0\leq t\leq T_1\eqdef c_2^{-1}<1$, 
\be\ba
\n{B(t)}{W^{1,q}} &\leq Cc_0\Big( 1+\int_0^t \| (B\cdot\g) v - B\div u\|_{W^{1,q}}ds\Big) \exp\Big(C\int_0^t\|\g v\|_{W^{1,q}} ds\Big)\\
& \leq Cc_0\Big( 1+\int_0^t \|B\|_{W^{1,q}}\|v\|_{W^{2,q}}ds\Big) \exp\Big(C\int_0^t\|v\|_{W^{2,q}} ds\Big)\\
& \leq Cc_0\Big( 1+\int_0^t \|B\|_{W^{1,q}}\|v\|_{W^{2,q}}ds\Big) \exp\Big(Cc_1^{\f 12}t^{\f 12}\Big)\\
& \leq Cc_0\Big( 1+\int_0^t \|B\|_{W^{1,q}}\|v\|_{W^{2,q}}ds\Big)
\ea\ee
which implies by the Gronwall inequality that
\[
\|B(t)\|_{W^{1,q}}\leq Cc_0^2\exp\Big( 1+\int_0^t\|v\|_{W^{2,q}}ds\Big)\leq Cc_0^2.
\]
Hence
\[
\|B_t(t)\|_{L^q}\leq C\|B(t)\|_{W^{1,q}}\|v(t)\|_{W^{1,q}} \leq Cc_0^2\|v(t)\|_{H^2} \leq Cc_0^2\beta c_1 = Cc_0^2c_2.
\]
\endproof

Finally, under the condition of $\rho_0\geq \delta >0$, one can prove the existence of $u$ to the linearized system \eqref{leq1}.
\begin{lemma}\label{LU}
Assume that $\rho_0\geq\delta>0$. Then, there exists a unique solution $u$ to \eqref{leq1} such that
\be\label{eu}
\left\{\ba
&\|u(t)\|_{H^1}^2 + \|\sqrt{\rho}u_t(t)\|_{L^2}^2 +\int_0^t\|u_t(s)\|_{H^1}^2 ds \leq Cc_0^{2\gamma+4}.\\
&\|\g u(t)\|_{H^1} \leq Cc_0^{\gamma+4}c_1^2,\\
&\int_0^t\|u(s)\|_{W^{2,q}}^2 ds\leq Cc_0^{2\gamma+8},
\ea\right.\ee
for $0\leq t\leq \min(T^*,T_2)$ with $T_2\eqdef \min(T_1, c_1^{-5}\exp(-Cc_2^{2\gamma+17}))$.
\end{lemma}
\proof
Direct energy equalities yield
\[\ba
\f 12\f{d}{dt}\int \rho u_t^2 dx& + \int \mu|\g u_t|^2 + (\mu + \ld)|\div u_t|^2 dx \\
&= \int (-\g p_t -\rho_t v\g u -\rho(2v\cdot\g u_t + v_t\cdot \g u) + \div(BB_t I - 2B\otimes B_t)) u_t dx
\ea\]
It follows from the  H\"older, Sobolev and Ladyzhenskaya inequalities that
\be\label{eu1}
\ba
&\f{d}{dt}\int \rho u_t^2 dx + \mu\int|\g u_t|^2dx \\
&\leq C\Big(\|p_t\|_{L^2}^2 + \|\rho_t\|_{L^2}^2 \|v\|_{L^\infty}^2 \|\g u\|_{L^4}^2 
	+ \|\sqrt{\rho}\|_{L^\infty}^2 \|v\|_{L^\infty}^2 \|\sqrt{\rho} u_t\|_{L^2}^2 +\|B\|_{L^\infty}\|B_t\|_{L^2}\Big) \\
&\qquad+ \|\sqrt{\rho}\|_{L^\infty} \|v_t\|_{L^4} \|\g u\|_{L^4}  \|\sqrt{\rho} u_t\|_{L^2}\\
&\leq C\Big(\|p_t\|_{L^2}^2 + \|\rho_t\|_{L^2}^2 \|v\|_{L^\infty}^2 \|\g u\|_{L^4}^2
	+ \|\sqrt{\rho}\|_{L^\infty}^2 \|v\|_{L^\infty}^2 \|\sqrt{\rho} u_t\|_{L^2}^2+\|B\|_{L^\infty}\|B_t\|_{L^2}\Big) \\
&\qquad + \f{C}{\eps}\|\rho\|_{L^\infty} \|\g u\|_{L^2} \|\g u\|_{H^1} +\eps \|\g v_t\|_{L^2}^2 \|\sqrt{\rho} u_t\|_{L^2}^2
\ea\ee
for any $\eps>0$. Noting that

\be\label{eu2}
\f 12 \f{d}{dt}\int |\g u|^2 dx = \int \g u\cdot \g u_t dx \leq \f{\mu}{2}\int |\g u_t|^2 dx + \f{1}{2\mu}\int |\g u|^2dx
\ee
and using the elliptic regularity theory, one can get
\be\label{eu3}
\ba
\|\g u\|_{H^1} &\leq C\left( \|\rho u_t\|_{L^2} + \|\rho v\cdot \g u\|_{L^2} + \|\g p\|_{L^2} + \|B\|_{W^{1,q}}\|\g B\|_{L^2}\right).\\
&\leq C\left(c_0^{\f 12}\|\sqrt{\rho} u_t\|_{L^2} + c_0c_2\|\g u\|_{L^2} + c_0^{\gamma+5}\right)
\ea\ee
Then
\[\ba
\|\rho_t\|_{L^2}^2 &\|v\|_{L^\infty}^2 \|\g u\|_{L^4}^2 
 	\leq C\left(c_2^4c_2^2 \|\g u\|_{H^1}^2\right)= Cc_2^6 \|\g u\|_{H^1}^2\\
&\leq Cc_2^6 \left(c_0^{\f 12}\|\sqrt{\rho} u_t\|_{L^2} + c_0c_2\|\g u\|_{L^2} + c_0^{\gamma+4}\right)^2\\
&\leq Cc_2^6 \left(c_0\|\sqrt{\rho} u_t\|_{L^2}^2 + c_0^2c_2^2\|\g u\|_{L^2}^2 + c_0^{2\gamma+8}\right)\\
&\leq Cc_0c_2^6 \|\sqrt{\rho} u_t\|_{L^2}^2 + Cc_0^2c_2^8\|\g u\|_{L^2}^2 + Cc_0^{2\gamma+8}c_2^6.
\ea\]

Combing \eqref{eu1}, \eqref{eu2} and \eqref{eu3} yields
\be\ba
\f{d}{dt}&\int \rho u_t^2 + |\g u|^2dx + \mu\int|\g u_t|^2dx \\
&\leq C\gamma^2c_0^{2\gamma-2}c_2^4 + Cc_2^6 \|\g u\|_{H^1}^2 + Cc_0c_2^2\|\sqrt{\rho}u_t\|_{L^2}^2 + Cc_0^4c_2\\
&\qquad+\f{C}{\eps}c_0\|\g u\|_{L^2}\left(c_0^{\f 12}\|\sqrt{\rho}u_t\|_{L^2} + c_0c_2\|\g u\|_{L^2} + c_0^{\gamma+4}\right)\\
&\qquad+ \eps\|\g v_t\|_{L^2}^2 \|\sqrt{\rho} u_t\|_{L^2}^2\\
&\leq C\gamma^2c_0^{2\gamma-2}c_2^4 + Cc_0c_2^6\|\sqrt{\rho}u_t\|_{L^2}^2 + Cc_0^2c_2^8\|\g u\|_{L^2}^2 + C c_0^{2\gamma+8}c_2^6 \\
	&\qquad+\f{C}{\eps}c_0\left(\|\g u\|_{L^2}^2 + c_0\|\sqrt{\rho}u_t\|_{L^2}^2 + c_0^2c_2^2\|\g u\|_{L^2}^2 + c_0^{2\gamma+8}\right)\\
	&\qquad+ \eps\|\g v_t\|_{L^2}^2 \|\sqrt{\rho} u_t\|_{L^2}^2\\
&\leq Cc_0^{2\gamma}c_2^4 + Cc_0c_2^6\|\sqrt{\rho}u_t\|_{L^2}^2 + Cc_0^2c_2^8\|\g u\|_{L^2}^2 + C c_0^{2\gamma+8}c_2^6 \\
	&\qquad+\f{C}{\eps}c_0\left(\|\g u\|_{L^2}^2 + c_0\|\sqrt{\rho}u_t\|_{L^2}^2 + c_0^2c_2^2\|\g u\|_{L^2}^2 + c_0^{2\gamma+8}\right)\\
	&\qquad+ \eps\|\g v_t\|_{L^2}^2 \|\sqrt{\rho} u_t\|_{L^2}^2\\
&\leq C\left(c_0^{2\gamma+8}c_2^6 + \f{c_0}{\eps}c_0^{2\gamma+8}\right) + C(c_0c_2^6+\f{c_0^2}{\eps})\|\sqrt{\rho}u_t\|_{L^2}^2 + C(c_0^2c_2^8+\f{c_0^3c_2^2}{\eps})\|\g u\|_{L^2}^2  \\
	&\qquad+ \eps\|\g v_t\|_{L^2}^2 \|\sqrt{\rho} u_t\|_{L^2}^2\\
&\leq \f{C}{\eps}c_2^{2\gamma+16}\left(1 + \|\sqrt{\rho}u_t\|_{L^2}^2 +\|\g u\|_{L^2}^2\right) + \eps\|\g v_t\|_{L^2}^2 \|\sqrt{\rho} u_t\|_{L^2}^2.
\ea\ee
 
 Taking $\eps = c_1^{-1}$, $t\leq T_2=\min(T_1, c_1^{-5}\exp(-Cc_2^{2\gamma+17}))$, and 
 \be\label{ERU1}
 \lim\limits_{t\ra 0+}\|\sqrt{\rho}u_t\|^2\leq Cc_0^{5},
 \ee
One gets by Gronwall's inequality the first estimate in \eqref{eu}.

To prove \eqref{ERU1}, one can rewrite the second equation of \eqref{leq1} as
\[
\rho u_t + \rho v\cdot \g u + \g p = \mu\div \D(u) + \ld\g\div u - B\times(\g\times B).
\]
That is
\[\ba
\rho^{\f 12} u_t&= - \rho^{\f 12} v\cdot \g u - \rho^{\f 12}\g p  +\mu\div \D(u) + \ld\g\div u -\rho^{\f 12} B\times(\g\times B).
\ea\]
Therefore,  it holds that
\[\ba
\lim\limits_{t\ra 0+}\|\rho^{\f 12} u_t\|_{L^2}
&=\lim\limits_{t\ra 0+}\|- \rho^{\f 12} v\cdot \g u - \rho^{\f 12}\g p  +\mu\div \D(u) + \ld\g\div u -\rho^{\f 12} B\times(\g\times B)\|_{L^2}\\
&=\|- \rho_0^{\f 12} u_0\cdot \g u_0 + g\|_{L^2}\\
&\leq \|\rho_0\|_{L^\infty}^{\f 12}\|u_0\|_{L^\infty}\|\g u_0\|_{L^2} + \| g\|_{L^2}\\
&\leq Cc_0^{\f{5}{2}}
\ea\]

Now, we prove the second estimate in \eqref{eu}. Since
\[\ba
\|\g u\|_{H^1} &\leq C\left( \|\rho u_t\|_{L^2} + \|\rho v\cdot \g u\|_{L^2} + \|\g p\|_{L^2} + \|B\|_{W^{1,q}}\|\g B\|_{L^2}\right).\\
&\leq C\left(c_0^{\f 12}\|\sqrt{\rho} u_t\|_{L^2} + c_0\|v\|_{L^4}\|\g u\|_{L^4} + c_0^{\gamma+4}\right)\\
&\leq C\left(c_0^{\gamma+4} + c_0c_1\|\g u\|_{L^2}^{\f 12} \|\g u\|_{H^1}^{\f 12}\right),
\ea\]
we have
\[\ba
\|\g u\|_{H^1} 
&\leq C\left(c_0^{\gamma+4} + c_0^2c_1^2\|\g u\|_{L^2}\right)\\
&\leq Cc_0^{\gamma+4}c_1^2.
\ea\]
which implies the second inequality in \eqref{eu}.

Finally, the standard elliptic regularity theory implies
\[\ba
\|u\|_{W^{2,q}}&\leq C\left(\|\rho u_t\|_{L^q} + \|\rho v\cdot \g u\|_{L^q} + \|\g p\|_{L^q} + \|B\times(\g\times B)\|_{L^q}\right)\\
&\leq C\left(c_0\|\g u_t\|_{L^2} + c_0c_2\|\g u\|_{L^q} + c_0^{\gamma+4}\right)\\
&\leq C\left(c_0\|\g u_t\|_{L^2} + c_0^{\gamma+4}\right).
\ea\]
Then
\[\ba
\int_0^{T_2}\|u\|_{W^{2,q}}^2ds&\leq C\left(c_0c_0^{\gamma+4}c_1^2 + c_0^{\gamma+4}\right)^2 T_2.
\ea\]
here $T_2=\min(T_1, c_1^{-5}\exp(-Cc_2^{2\gamma+17}))$.
This, together with $\eqref{eu}_1$, gives $\eqref{eu}_3$.
\endproof

Now, let $c_1\eqdef 4Cc_0^{2\gamma+8},\ c_2\eqdef\f 12Cc_0^{\gamma+6}c_1^2=8C^3c_0^{5\gamma+22}$ and $\beta\eqdef\f{c_2}{c_1}=4C^2c_0^{3\gamma+14}$, we have
\be\label{AE1}
\left\{\ba
& \|\rho\|_{W^{1,q}} + \|\rho_t\|_{L^q} + \|B\|_{W^{1,q}} + \|B_t\|_{L^q} \leq Cc_2^{2\gamma+4},\\
& \|\sqrt{\rho}u_t(t)\|_{L^2} \leq Cc_0^{2\gamma+4}\leq Cc_2^{2\gamma+4},\\
& \|u(t)\|_{H^1} + \beta^{-1}\|u(t)\|_{H^2} + \int_0^t \|u_t(s)\|_{H^1}^2 + \|u(s)\|_{W^{2,q}}^2 ds \leq c_1.
\ea\right.\ee
for $0\leq t\leq \min(T^*, T_2)$.

Using the same arguments in \cite{2006ChoK}, we obtain
\begin{lemma}\label{LP1} There exists a unique solution $(\rho, u, B)$ to the linear problem \eqref{leq1} in $[0, T_*]$ satisfying  \eqref{AE1} as well as the following regularity
\be\label{LR1}\left\{\ba
& \rho \in C([0, T_*]; W^{1,q}),\ \rho_t\in C([0, T_*]; L^q),\\
& (u, B)\in C([0, T_*]; H^2)\cap L^2(0, T_*; W^{2,q}),\\
& (u_t,B_t)\in L^2(0, T_*; H^1),\\
& (\sqrt{\rho} u_t, B_t)\in L^{\infty}(0, T_*; L^2).
\ea\right.\ee
where $T_*=\min(T^*, T_2)$.
\end{lemma}
\proof For any $\delta\in (0,1)$, set $\rho_0^\delta = \rho_0+\delta$, $B_0^\delta=B_0$. Let $u_0^\delta\in H_0^1\cap H^2$ be the solution of
\[\left\{\ba
&\div(\mu\D( u_0^\delta)) + \g(\ld \div u_0^\delta) = \g p(\rho_0^\delta)  + B_0^\delta\times(\g\times B_0^\delta) - (\rho_0^\delta)^{\f 12}g,\\
&u_0^\delta = 0,\quad \text{ on } \p\O.
\ea\right.\]
Then it holds that
\[
\g p(\rho_0^{\delta}) - \div(\mu\D( u_0^\delta)) - \g(\ld \div u_0^\delta) + B_0^\delta\times(\g\times B_0^\delta) = (\rho_0^\delta)^{\f 12}g
\]
and
\[
\|u_0^\delta - u_0\|_{H^2}\leq C\left(\|p(\rho_0^\delta)-p(\rho_0)\|_{H^1} + \|((\rho_0^\delta)^{\f 12}-\rho_0^{\f 12})g\|_{L^2}\right)\ra 0.
\]
Moreover, for all $\delta>0$ small enough, we have
\[
2+\gamma + \delta + \|\rho_0^\delta-\delta\|_{W^{1,q}} + \n{u_0^\delta}{H^2} +\n{g}{L^2}^2 < c_0.
\]

Then, it follows from Lemmas \ref{LRHO}-\ref{LU} that, for any small $\delta>0$, there exists unique strong solution $(\rho^\delta, u^\delta, B^\delta)$ to the  linearized problem \eqref{leq1}-\eqref{leq3} which satisfies the estimates \eqref{AE1}. It follows from standard compactness arguments that, there exists a subsequence of solutions $(\rho^\delta, u^\delta, B^\delta)$ converges to $(\rho, u, B)\in C([0, T], H^1)$. 

Clearly, $(\rho, u, B)\in C([0, T_2], H^1)$ also satisfies the estimates in \eqref{AE1}.

It is easy to show that $(\rho, u, B)$ is a strong solution to the linearized problem \eqref{leq1}-\eqref{leq3} in $[0, T_*)$ satisfying the  regularity \eqref{LR1}.

Now, we prove the uniqueness of the solution in this regularity class. Let $(\rho_1, u_1, B_1)$ and $(\rho_2, u_2, B_2)$ be two such solutions. Set
\[
\t\rho = \rho_1 - \rho_2,\quad \t\u = u_1-u_2,\quad \t B=B_1-B_2.
\]

Then, $\t\rho\in L^\infty(0,T_*; L^q)$ solves the following transport equation 
\[
\t\rho_t + \div(\t\rho v) = 0.
\]
Therefore $\t\rho=0$. That is $\rho_1=\rho_2$. Similiarly $B_1=B_2$.

Multiplying the both sides of 
\[
\rho\t u_t + \rho v\cdot \t u = \div(\mu\D(\t u)) + \g(\ld\div \t u)
\]
by $\t u$ and integrating by part, one has
\[
\f 12\int\rho |\t u|^2dx + \int_0^t\int \mu|\D(\t u)|^2 + \ld|\div \t u|^2 dxds = 0.
\]
This gives $u_1=u_2$.
\endproof

\section{Nonlinear Problem}
Base on some a priori estimates, one can prove the existence and regularity of a unique solution $(\rho, u, B)$ to the nonlinear problem as follows.

{\bf Proof of theorem 1.} Let $u^{(0)}\in C([0,\infty; H^2\cap H_0^1)$ be the solution to
\[
w_t -\lap w = 0,\quad w(0)= u_0,\ w|_{\p\O}=0.
\]
Then
\[
\sup\limits_{0\leq t\leq T_*}\left(\|u^{(0)}(t)\|_{ H_0^1} + \beta^{-1}\|u^{(0)}(t)\|_{H^2}\right) + \int_0^{T_*} \|u^{(0)}_t(s)\|_{H^1}^2 + \|u^{(0)}(s)\|_{W^{2,q}}^2ds \leq c_1.
\]
By Lemma \ref{LP1}, there exists a unique strong solution $(\rho^{(1)}, u^{(1)}, B^{(1)})$ to the linearized problem \eqref{leq1}-\eqref{leq3} with
$v=u^{(0)}$. Then, one can construct approximate solutions $(\rho^{(i)}, u^{(i)}, B^{(i)})$ by induction. Assume that $u^{(i-1)}$ have been defined for $i\geq 1$. Let $(\rho^{(i)}, u^{(i)}, B^{(i)})$ be the 
unique solution to the linearized problem \eqref{leq1}-\eqref{leq3} with $v=u^{(i-1)}$.

Then, it follows from Lemma \ref{LP1} that there exists a constant $C>1$ such that
\[\left\{\ba
& \|\rho^{(i)}\|_{W^{1,q}} + \|\rho^{(i)}_t\|_{L^q} + \|B^{(i)}\|_{W^{1,q}} + \|B^{(i)}_t\|_{L^q} \leq Cc_2^{2\gamma+4},\\
& \|\sqrt{\rho^{(i)}}u^{(i)}_t(t)\|_{L^2} \leq Cc_0^{2\gamma+4}\leq Cc_2^{2\gamma+4},\\
& \|u^{(i)}(t)\|_{H^1} + \beta^{-1}\|u^{(i)}(t)\|_{H^2} + \int_0^t \|u^{(i)}_t(s)\|_{H^1}^2 + \|u^{(i)}(s)\|_{W^{2,q}}^2 ds \leq c_1.
\ea\right.\]

Now, we will prove that $(\rho^{(i)}, u^{(i)}, B^{(i)})$ converges to a solution to the nonlinear problem in a strong sense.

Let
\[
\t\rho^{(i+1)} = \rho^{(i+1)}-\rho^{(i)}, \t u^{(i+1)} =  u^{(i+1)} - u^{(i)}, \t B^{(i+1)} =  B^{(i+1)} - B^{(i)}.
\]
Then
\be\label{C1}\left\{\ba
&\t\rho^{(i+1)}_t + \div(\t\rho^{(i+1)} u^{(i)}) + \div(\rho^{(i)}\t u^{(i)}) = 0,\\
&\rho^{(i+1)}\t u^{(i+1)}_t + \rho^{(i+1)} u^{(i)}\cdot\g \t u^{(i+1)} -\mu\div\D(\t u^{(i+1)}) -\ld\g\div u^{(i+1)} \\
	&\qquad =\t\rho^{(i+1)}(-u^{(i)}_t-u^{(i-1)}\cdot\g u^{(i)}) -\rho^{(i+1)}\t u^{(i)}\cdot\g u^{(i)} -\g(p(\rho^{(i+1)})-p(\rho^{(i)}))\\
	&\qquad\qquad - B^{(i+1)}\times(\g\times \t B^{(i+1)}) - \t B^{(i+1)}\times(\g\times B^{(i)}) =0,\\
&\t B^{(i+1)}_t +\g\times (\t B^{(i+1)}\times  u^{(i)}+ B^{(i)}\times \t u^{(i)})=0,\\
&\div \t B^{(i+1)} = 0.
\ea\right.\ee

Multiplying $\eqref{C1}_1$ by $\t\rho^{(i+1)}$ and integrating by part lead to
\[\ba
\f{d}{dt}\int |\t\rho^{(i+1)}|^2 dx &\leq C\int |\g u^{(i)}||\t\rho^{(i+1)}|^2 + |\g\rho^{(i)}| |\t u^{(i)}| |\t\rho^{(i+1)}| + |\rho^{(i)}| |\g\t u^{(i)}| |\t\rho^{(i+1)}| dx\\
&\leq C\left(\|\g u^{(i)}\|_{W^{1,q}}\|\t\rho^{(i+1)}\|_{L^2}^2 +(\|\g\rho^{(i)}\|_{L^q} +\|\rho^{(i)}\|_{L^\infty})\|\g\t u^{(i)}\|_{L^2} \|\t\rho^{(i+1)}\|_{L^2}\right).
\ea\]
Hence, by Young's inequality, one gets
\be\label{C1}
\f{d}{dt}\int |\t\rho^{(i+1)}|^2 dx \leq A_\eps^{(i)} \|\t\rho^{(i+1)}\|_{L^2}^2 +\eps\|\g\t u^{(i)}\|_{L^2}^2
\ee
where
\[
A_\eps^{(i)} \eqdef C\|\g u^{(i)}\|_{W^{1,q}} + \f{C}{\eps}(\|\g\rho^{(i)}\|_{L^q}^2 +\|\rho^{(i)}\|_{L^\infty}^2).
\]

Multiplying $\eqref{C1}_3$ by $\t B^{(i+1)}$ and integrating by part yield
\[\ba
\f{d}{dt}\|\t B^{(i+1)}\|_{L^2}^2 & \leq C\|\g u^{(i)}\|_{L^\infty} \|\t B^{(i+1)}\|_{L^2}^2 \\
&\qquad +	C\left(\|\g B^{(i)}\|_{L^q} + \|B^{(i)}\|_{L^\infty}\right) \|\g\t u^{(i)}\|_{L^2}\|\t B^{(i+1)}\|_{L^2}.
\ea\]
Hence, by Young's inequality one gets
\be\label{C2}
\f{d}{dt}\|\t B^{(i+1)}\|_{L^2}^2 \leq B_\eps^{(i)}\|\t B^{(i+1)}\|_{L^2}^2 +\eps  \|\g\t u^{(i)}\|_{L^2}^2,
\ee
where
\[
B_\eps^{(i)} \eqdef C\|\g u^{(i)}\|_{L^\infty} +\f{C}{\eps}\left(\|\g B^{(i)}\|_{L^q}^2 + \|B^{(i)}\|_{L^\infty}^2\right).
\]

Multiplying $\eqref{C1}_2$ by $\t u^{(i+1)}$   and integrating by part, one can get
\[\ba
\f12 \f{d}{dt}&\int\rho^{(i+1)}|\t u^{(i+1)}|^2 dx + \mu\int |\D(\t u^{(i+1)})|^2 + \ld|\div \t u^{(i+1)}|^2dx\\
&\leq C\int |\t\rho^{(i+1)}| (|u^{(i)}_t| + |u^{(i-1)}\g u^{(i)}|) |\t u^{(i+1)}| + |\rho^{(i+1)}| |\t u^{(i)}| |\g u^{(i)}| |\t u^{(i+1)}|\\
	&\qquad + |p(\rho^{(i+1)})-p(\rho^{(i)})| |\g\t u^{(i+1)}| + (|B^{(i)}| + |B^{(i+1)}|) |\t B^{(i+1)}| |\g\t u^{(i+1)}|dx
\ea\]
Then
\[\ba
\f 12\f{d}{dt}&\|\sqrt{\rho^{(i+1)}}\t u^{(i+1)}\|_{L^2}^2 + \mu\|\g\t u^{(i+1)}\|_{L^2}^2\\
&\leq C\|\t\rho^{(i+1)}\|_{L^2} (\|u^{(i)}_t\|_{L^4} + \|u^{(i-1)}\|_{L^\infty}\|\g u^{(i)}\|_{L^4}) \|\t u^{(i+1)}\|_{L^4} \\
	&\qquad + \|\rho^{(i+1)}\|_{L^\infty}^{\f 12} \|\t u^{(i)}\|_{L^4} \|\g u^{(i)}\|_{L^4} \|\sqrt{\rho^{(i+1)}}\t u^{(i+1)}\|_{L^2}\\
	&\qquad + \|p(\rho^{(i+1)})-p(\rho^{(i)})\|_{L^2} \|\g\t u^{(i+1)}\|_{L^2}\\
	&\qquad + (\|B^{(i)}\|_{L^\infty} + \|B^{(i+1)}\|_{L^\infty}) \|\t B^{(i)}\|_{L^2} \|\g\t u^{(i+1)}\|_{L^2}\\
&\leq C\|\t\rho^{(i+1)}\|_{L^2} \|\g\t u^{(i+1)}\|_{L^2} (\|u^{(i)}_t\|_{L^4} + 1) \\
	&\qquad + C \|\g\t u^{(i)}\|_{L^2} \|\sqrt{\rho^{(i+1)}}\t u^{(i+1)}\|_{L^2}\\
	&\qquad + C\|\t\rho^{(i+1)}\|_{L^2} \|\g\t u^{(i+1)}\|_{L^2}\\
	&\qquad + C\|\t B^{(i+1)}\|_{L^2} \|\g\t u^{(i+1)}\|_{L^2}.
\ea\]
It then follows from this and Young's inequality that
\[\ba
\f 12\f{d}{dt}&\|\sqrt{\rho^{(i+1)}}\t u^{(i+1)}\|_{L^2}^2 + \mu\|\g\t u^{(i+1)}\|_{L^2}^2\\
&\leq \f{C}{\mu}\|\t\rho^{(i+1)}\|_{L^2}^2  (\|u^{(i)}_t\|_{L^4} + 1)^2 \\
	&\qquad + \f{C}{\eps} \|\sqrt{\rho^{(i+1)}}\t u^{(i+1)}\|_{L^2}^2\\
	&\qquad + \f{C}{\mu}\|\t\rho^{(i+1)}\|_{L^2}^2 \\
	&\qquad + \f{C}{\mu}\|\t B^{(i+1)}\|_{L^2}^2\\
	& +\eps \|\g\t u^{(i+1)}\|_{L^2}^2 + \f{\mu}{2}\|\g\t u^{(i)}\|_{L^2}^2.
\ea\]
Hence, 
\be\label{C3}\ba
\f 12\f{d}{dt}&\|\sqrt{\rho^{(i+1)}}\t u^{(i+1)}\|_{L^2}^2 + \f{\mu}{2}\|\g\t u^{(i+1)}\|_{L^2}^2\\
&\leq C\|\t\rho^{(i+1)}\|_{L^2}^2  (\|u^{(i)}_t\|_{L^4} + 1)^2
	 + \f{C}{\eps} \|\sqrt{\rho^{(i+1)}}\t u^{(i+1)}\|_{L^2}^2\\
	 &\qquad + C\|\t B^{(i+1)}\|_{L^2}^2 +\eps\|\g\t u^{(i)}\|_{L^2}^2.
\ea\ee

Summarizing \eqref{C1}-\eqref{C3}, one has obtained that
\[\left\{\ba
\f{d}{dt}&\|\t\rho^{(i+1)}\|_{L^2}^2 \leq A_\eps^{(i)} \|\t\rho^{(i+1)}\|_{L^2}^2 +\eps\|\g\t u^{(i)}\|_{L^2}^2\\
\f{d}{dt}&\|\t B^{(i+1)}\|_{L^2}^2 \leq B_\eps^{(i)}\|\t B^{(i+1)}\|_{L^2}^2 +\eps  \|\g\t u^{(i)}\|_{L^2}^2,\\
\f{d}{dt}&\|\sqrt{\rho^{(i+1)}}\t u^{(i+1)}\|_{L^2}^2 + \mu\|\g\t u^{(i+1)}\|_{L^2}^2\\
&\leq C\|\t\rho^{(i+1)}\|_{L^2}^2  (\|u^{(i)}_t\|_{L^4} + 1)^2 
+ \f{C}{\eps} \|\sqrt{\rho^{(i+1)}}\t u^{(i+1)}\|_{L^2}^2
+ C\|\t B^{(i+1)}\|_{L^2}^2 +\eps\|\g\t u^{(i)}\|_{L^2}^2.
\ea\right.\]
with
\[\ba
&A_\eps^{(i)} \eqdef C\|\g u^{(i)}\|_{W^{1,q}} + \f{C}{\eps}(\|\g\rho^{(i)}\|_{L^q}^2 +\|\rho^{(i)}\|_{L^\infty}^2),\\
&B_\eps^{(i)} \eqdef C\|\g u^{(i)}\|_{L^\infty} +\f{C}{\eps}\left(\|\g B^{(i)}\|_{L^q}^2 + \|B^{(i)}\|_{L^\infty}^2\right).
\ea\]

Define 
\[
\psi^{(i+1)}(t)\eqdef \|\t\rho^{(i+1)}\|_{L^2}^2 + \|\t B^{(i+1)}\|_{L^2}^2  + \|\sqrt{\rho^{(i+1)}}\t u^{(i+1)}\|_{L^2}^2.
\]
Then
\[\ba
\f{d}{dt}&\psi^{(i+1)}(t)  + \mu\|\g\t u^{(i+1)}\|_{L^2}^2\\
	&\leq C\left( A_\eps^{(i)} +  B_\eps^{(i)} +\f{C}{\eps} + \|u^{(i)}_t\|_{L^4}^2 + 1  \right)\psi^{(i+1)} + \eps\|\g\t u^{(i)}\|_{L^2}^2.
\ea\]
Now, recalling that $\psi^{(i+1)}(0)=0$ and using Gronwall's inequality, we have
\[\ba
\psi^{(i+1)}(t) &+ \mu\|\g\t u^{(i+1)}\|_{L^2}^2 ds\\
&\leq \eps\left(1+\int_0^t F^{(i)}_\eps(s)ds\exp\left(\int_0^t F^{(i)}_\eps(s)ds\right)\right)\int_0^t\|\g\t u^{(i)}\|_{L^2}^2ds
\ea\]
with
\[
F^{(i)}_\eps(t) = C\left( A_\eps^{(i)} +  B_\eps^{(i)} +\f{C}{\eps} + \|u^{(i)}_t\|_{L^4}^2 + 1  \right).
\]

Since
\[
\int_0^t F^{(i)}_\eps(s)ds \leq C t + C\sqrt{t} + \f{C}{\eps}t + C\eps t + C\eps + C,
\]
one can choose $T_{**}\leq \eps<1$ so that
\[\ba
\psi^{(i+1)}(t) & + \mu\int_0^t\|\g\t u^{(i+1)}\|_{L^2}^2 ds\\
&\leq \exp(C)\eps\int_0^t\|\g\t u^{(i)}\|_{L^2}^2ds
\ea\]
for $0\leq t\leq T^*\eqdef \min(T_*, T_{**})$. 

Then, taking $\eps$ suitably small so that 
\[
\eps\exp(C)\leq \f1 2\mu,
\]
one obtains
\[\ba
\sum\limits_{i=0}^{\infty}\sup\limits_{0\leq t\leq  T^*}\psi^{(i+1)}(t) &+ \sum\limits_{i=0}^{\infty}\int_0^{ T^*} \mu\|\g\t u^{(i+1)}\|_{L^2}^2 ds\\
&\leq \sum\limits_{i=0}^{\infty}\exp(C)\eps\int_0^{ T^*}\|\g\t u^{(i)}\|_{L^2}^2ds
\ea\]
Hence
\[\ba
\sum\limits_{i=0}^{\infty}\sup\limits_{0\leq t\leq  T^*}\psi^{(i+1)}(t) &+ \sum\limits_{i=0}^{\infty}\int_0^{ T^*} \mu\|\g\t u^{(i+1)}\|_{L^2}^2 ds\\
&\leq \mu\int_0^{ T^*}\|\g\t u^{(0)}\|_{L^2}^2ds\leq C.
\ea\]
Therefore, $(\rho^{(i)}, u^{(i)}, B^{(i)})$ converges to a limit $(\rho, u, B)$ in the following sense:
\[\ba
&(\rho^{(i)}, B^{(i)}) \ra (\rho, B), \text{ in } L^{\infty}(0, T^*; L^2),\\
&u^{(i)} \ra u \text{ in } L^2(0, T^*; H^1).
\ea\]

It is easy to prove that $(\rho, u, B)$ is a weak solution to the nonlinear problem and satisfies the following estimate
\[\ba
\sup\limits_{0\leq t\leq  T^*}&\|\sqrt{\rho} u_t\|_{L^2} + \int_0^{ T^*}\|u_t\|_{L^2}^2 + \|u_t\|_{W^{2,q}}^2 dt \\
&\qquad +\sup\limits_{0\leq t\leq  T^*}\left(\|(\rho, B)\|_{W^{1,q}} + \|(\rho_t, B_t)\|_{L^q} + \|u\|_{H^2}\right) \leq C.
\ea\]
This proves the local existence of a strong solution. The uniqueness of the strong solutions is trivial thus omitted here.
\endproof

\section*{Acknowledgements}
 X.-D. Huang is partially supported by CAS Project for Young Scientists in Basic Research, Grant No.YSBR-031 and NNSFC Grant No. 11688101. Xin is partially supported by Zheng Ge Ru Foundation, Hong Kong RGC Earmarked Research Grants CUHK-14301421, CUHK-14301023, CUHK-14300819 and CUHK-14302819, and the key projects of NSFC Gronts No. 12131010 and No. 11931013. W. YAN is partially supported by NNSFC Grant No. 11371064 and No. 11871113.

\end{document}